\documentclass[12pt,reqno,a4wide]{amsart}

\usepackage{tikz} 
\usepackage{calrsfs}
\usepackage[charter]{mathdesign}
\usepackage{hyperref}

\usetikzlibrary{decorations.pathreplacing}
\usetikzlibrary{fit}

\usepackage{pgfplots}

\allowdisplaybreaks

\newtheorem*{defi}{Definition}

\newtheorem*{acka}{Acknowledgment}

\catcode`,\active

\catcode`\,12

\begin{document}
\bibliographystyle{plain}

%
%

	\title[Enumeration   of  S-Motzkin paths]
	{An elementary approach to solve recursions relative to the enumeration   of  S-Motzkin paths }

	\author[H. Prodinger ]{Helmut Prodinger }
	\address{Department of Mathematics, University of Stellenbosch 7602, Stellenbosch, South Africa}
	\email{hproding@sun.ac.za}

	
	\begin{abstract}
S-Motzkin paths (bijective to ternary trees) and partial version of them are calculated using only
elementary methods from linear algebra.
	\end{abstract}
	

\maketitle


%
%

\section{Introduction}

\emph{Motzkin} paths consist, as their cousins, the Dyck paths, of up steps and down steps, but additionally of horizontal steps. They start at the origin, end at the $x$-axis, and never go below the $x$-axis. The encyclopedia \cite{OEIS} contains more than 600 related items.

The subfamily of Motzkin paths called  S-Motzkin paths has been introduded in \cite{Selkirk-arxiv}  and analyzed in \cite{PSW}. Further papers on the subject are \cite{Hana, Nancy-bijection}. 

S-Motzkin paths are Motzkin paths that consist of $n$ up steps, $n$ horizontal steps, and $n$ down steps. Ignoring down steps ($d$), the sequence of horizontal ($h$) and up ($u$) steps must look like $huhu\dots hu$. They are enumerated by
\begin{equation*}
\frac1{2n+1}\binom{3n}{n}.
\end{equation*}

In \cite{alone} \emph{partial} S-Motzkin paths were introduced and analyzed: If one stops the left-to-right traversal of an S-Motzkin path at some time, we might still be at level $k$; since flat steps and up steps interchange, one needs to introduce two families to   analyze this more general scenario of not necessarily ending on the $x$-axis, and also one can consider the S-Motzkin paths from right to left, ending at any prescribed level.

Let $a_{n, k}$ denote the number of partial S-Motzkin paths of length $n$ that  end at height $k$ and the last non-down step is an up step.  Furthermore, let $b_{n, k}$ denote the number of partial S-Motzkin path of length $n$ that end at height $k$ and the last  non-down step is a level  step.
The natural choice for the initial values are $a_{0, 0} = 1$ and $b_{0, 0} = 0$.

The analysis in \cite{alone} is appealing and instructive but uses heavy analytic machinery. The present paper aims at a (relatively) elementary treatment, using only systems of linear equations and Cramer's rule to solve them. Some determinants appear, and they follow some third-order recursions that can be solved. The help of \textsc{Maple} and the package \textsc{Gfun}~\cite{gfun} is gratefully acknowledged.

As one reviewer suggested, we briefly describe the approach chosen in the earlier paper \cite{alone}. At the core of the analysis are the bivariate generating functions
$A(z,u)=\sum_{n,k}a_{n,k}z^nu^k$ and $B(z,u)=\sum_{n,k}b_{n,k}z^nu^k$, which turn out to be  rational functions with a cubic denominator in the variable $u$. The \emph{kernel method} allows then to cancel one of the factors from the denominator and this makes the generating function quite explicit. Following the same approach with $C(z,u)=\sum_{n,k}c_{n,k}z^nu^k$ and $D(z,u)=\sum_{n,k}d_{n,k}z^nu^k$ (for related sequences $c_{n,k}$ and $d_{n,k}$, described briefly below and in more detail later), it comes as a surprise (at least in the beginning) that \emph{two} of three factors need to be cancelled. This leads in either case to a set of  two equations that can be solved.

It is easily seen that the following recurrence relations hold (details to follow):
\begin{align*}
a_{n, k} & = b_{n-1, k-1} + a_{n-1, k+1},\\
b_{n, k} & = a_{n-1, k} + b_{n-1, k+1}.
\end{align*} 
In the same spirit we get the other pair of recursions (details in a later section)
 \begin{align*}
c_{n, k} & = c_{n-1, k-1} + d_{n-1, k},\\
d_{n, k} & = d_{n-1, k-1} + c_{n-1, k+1}.
\end{align*}

We may note that $a_{n,0}=c_{n,0}$, since both just enumerate S-Motzkin paths; we get
\begin{equation*}
a_{3n,0}=c_{3n,0}=\frac1{2n+1}\binom{3n}{n}.
\end{equation*}

The following sections show how to solve these two systems of recursions by basically elementary methods.

By `solving' we mean to find explicit expressions for the respective generating functions. One can then go one step further
and read off coefficients, but that is not new, as it was already shown in \cite{alone}.

The following shows an S-Motzkin path of length $21$; the red part alone is a partial S-Motzkin path ending on level 3, 
and the blue (dashed) part alone, read from right to left is  a (reversed) partial S-Motzkin path ending on level 3. Section 2 is on the enumeration of the red objects, and Section 3 is on the enumeration of the blue (dashed) objects.

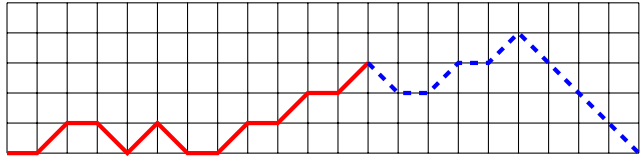
\begin{figure}[h]
\begin{center}
	\begin{tikzpicture}[scale=0.4]
	
	\draw[step=1.cm,black,thin ] (-0.0,-0.0) grid (21.0,5.0);

	\draw[ultra thick,red] (0,0) to (1,0) to (2,1) to (3,1) to (4,0) to (5,1) to (6,0) to (7,0) to (8,1) to (9,1)
	to (10,2)to (11,2)to (12,3);
	 
	 \draw[ultra thick,blue, dashed] (12,3)to (13,2)to (14,2)to (15,3)to (16,3)to (17,4)to (18,3)to (19,2)to (20,1) to (21,0);

	\end{tikzpicture}
\end{center}
\caption{An example of an S-Motzkin path of length 21, separated in a red part in the beginning and a blue (dashed) part at the end.}
\end{figure}

While the analysis in this paper is focused on the elementary analysis of partial S-Motzkin paths, we strongly believe that the approach with
determinants is of more general interest. Furthermore, in order to make this paper readable independently of \cite{alone}, a few repetitions were
unavoidable.

\section{The first pair of recursions}

The objects of this section are the numbers $a_{n,k}$, $b_{n,k}$ and their recurrences.
\begin{defi}
The numbers $a_{n,k}$\ ($b_{n,k}$) are defined as the number of partial S-Motzkin paths of length $n$, ending at level $k$ such the last non-down step
is an up (level) step.
\end{defi}
These sequences are coupled by these recursions:
\begin{align*}
	a_{n, k} & = b_{n-1, k-1} + a_{n-1, k+1},\\
	b_{n, k} & = a_{n-1, k} + b_{n-1, k+1}.
\end{align*} 
The recurrence relations are valid for $n\ge1$ and $k\ge0$; the quantity $b_{n-1, -1}$ must be interpreted as 0. Furthermore, the initial conditions that make sense
are $a_{0,0}=1$ and $b_{0,0}=0$. The recursions can be understood when distinguishing the cases of the last step leading to level $k$. Note further that $a_{n,0}$ is the number of S-Motzkin paths of length $n$.

We use generating functions
\begin{equation*}
f_k=f_k(z)=\sum_{n\ge0}a_{n,k}z^n\quad\text{and}\quad g_k=g_k(z)=\sum_{n\ge0}b_{n,k}z^n.
\end{equation*}
Then the sets of recursions translate into
\begin{align*}
f_{ k} & = zg_{ k-1} + zf_{ k+1}+[k=0],\\
g_{ k} & = zf_{ k} + zg_{ k+1},
\end{align*} 
where we employed Iverson's notation.\footnote{This convenient notation is quite common these days: $[P]=1$ if condition $P$ is true otherwise $[P]=0$.} Our goal is to eliminate one set of generating functions:
We find
\begin{align*}
f_{ k} -zf_{ k+1}& = -z(zg_{ k}-g_{ k-1}) + zf_{ k+1}-z^2f_{ k+2}+[k=0],\\
\end{align*}
and further
\begin{align*}
f_{ k}  & = z^2f_{k-1} + 2zf_{ k+1}-z^2f_{ k+2}+[k=0].
\end{align*}
These recursions are best expressed in matrix form. It requires an infinite matrix, but  we limit it by a parameter $h$, so that traditional methods from linear algebra can be used. In the following example $h=5$, since the functions $f_0,\dots,f_4$ are involved. At an appropriate time, we will push this parameter $h$ to infinity.
\begin{equation*}
\begin{pmatrix}
	1&-z&&&\\
	-z^2&1&-2z&z^2\\
&	-z^2&1&-2z&z^2\\
&&	-z^2&1&-2z&z^2\\
&&&	-z^2&1&-2z\\
\end{pmatrix}
\begin{pmatrix}
	f_0\\
		f_1\\
			f_2\\
				f_3\\
					f_4\\
\end{pmatrix}
=\begin{pmatrix}
	1\\
	0\\
	0\\
	0\\
	0\\
\end{pmatrix}
\end{equation*}
As we see, the first row  of the matrix is different. In order to overcome this, we introduce the following determinants (they all depend on the variable $z$, but we don't write this explicitly most of the time):
\begin{equation*}
D_h=D_h(z)=\det
\begin{pmatrix}
1&-2z&z^2&&\\
-z^2&1&-2z&z^2\\
&	-z^2&1&-2z&z^2\\
&&	-z^2&1&-2z&z^2\\
&&&	-z^2&1&-2z\\
\end{pmatrix}
\end{equation*}
and
\begin{equation*}
D_h^*=D_h^*(z)=\det
\begin{pmatrix}
1&-z&&&\\
-z^2&1&-2z&z^2\\
&	-z^2&1&-2z&z^2\\
&&	-z^2&1&-2z&z^2\\
&&&	-z^2&1&-2z\\
\end{pmatrix},
\end{equation*}
in both cases $h$ refers to the numbers of rows/columns. Expanding $D_h^*$ according to the first row, say, we find
\begin{equation*}
D_h^*=D_{h-1}-z^3D_{h-2},
\end{equation*}
and it is easier to work with the  set of determinants $D_k$, since the structure is more regular.
Expanding this determinant $D_k$ (again according to the first row, say) we get the recursion
\begin{equation*}
D_{k}-D_{k-1}+2z^3D_{k-2}-z^6D_{k-3}=0,
\end{equation*}
which has the characteristic equation $X^3-X^2+2z^3X-z^6=0$. In order to deal with the cubic equation successfully, it is 
a good idea to use an auxiliary variable: $z^3=t(1-t)^2$. This strategy has been applied before in other publications, for instance in
\cite{naiomi}. Then the 3 roots have a beautiful form:
\begin{align*}
r_1&=(t-1)^2,\quad
r_2=\frac t2(2-t+W),\quad
r_3=\frac t2(2-t-W),
\end{align*}
with $W=\sqrt{4t-3t^2}$. Consequently
 \begin{equation*}
D_h=Ar_1^h+Br_2^h+Cr_3^h,
\end{equation*}
and the constants will be worked out from the initial conditions $D_0=1$, $D_1=1$, $D_2=1-2z^3$. We find (of course using computer algebra!)
\begin{align*}
A&=\frac{t-1}{3t-1},\quad
B=\frac{3t^2-4t-W}{(3t-1)(3t-4)},\quad
C=\frac{3t^2-4t+W}{(3t-1)(3t-4)}.
\end{align*}
Now we will use these determinants to solve for the unknown generating functions $f_i$, using Cramer's rule. Allowing a finite determinant with $h$ rows and columns, means that we solve for the $h$ functions $f_0,f_1,\dots,f_{h-1}$, and we get 
 \begin{equation*}
f_i=\frac{z^{2i}D_{h-i-1}}{D_h^*}=\frac{z^{2i}D_{h-i-1}}{D_{h-1}-z^3D_{h-2}}.
 \end{equation*}
 Now it is time to consider the limit $h\to\infty$. It is not too difficult to see that the contributions from the roots $r_2$ and $r_3$ will disappear.
 This follows from the expansions 
 \begin{equation*}
r_{2,3}=-\frac12t^2\pm t^{3/2}\mp \frac38t^{5/2}+\cdots
 \end{equation*}
for small $t$, and $t\approx x$. In this way,  we find
 \begin{align*}
f_i&=\lim_{h\to\infty}\frac{z^{2i}D_{h-i-1}}{D_{h-1}-z^3D_{h-2}}
=\lim_{h\to\infty}\frac{z^{2i}Ar_1^{h-i-1}}{Ar_1^{h-1}-z^3Ar_1^{h-2}}\\
&=\frac{z^{2i}r_1^{-i-1}}{r_1^{-1}-z^3r_1^{-2}}=\frac{t^i}{z^i(1-t)}.
 \end{align*}
 From the system of recursions we further find
 \begin{align*}
zg_k=f_{k+1}-zf_{k+2}=\frac{t^{k+1}}{z^{k+1}(1-t)}-\frac{t^{k+2}}{z^{k+1}(1-t)}=\frac{t^{k+1}}{z^{k+1}},
 \end{align*}
 and finally 
 \begin{equation*}
g_k=\frac{t^{k+1}}{z^{k+2}}.
 \end{equation*}
 
 \section{The second pair of recursions}
A (partial) \textit{reverse S-Motzkin} path is a (partial) S-Motzkin path read from right to left.

\begin{defi}
	The numbers $c_{n,k}$\quad ($d_{n,k}$) are defined as the number of partial S-Motzkin paths of length $n$, ending at level $k$ such the last non-up step
	is an level (down) step.
\end{defi}

The system of recursions is
 \begin{align*}
	c_{n, k} & = c_{n-1, k-1} + d_{n-1, k},\\
	d_{n, k} & = d_{n-1, k-1} + c_{n-1, k+1};
\end{align*}
these recurrences hold for $n\ge1$ and $k\ge0$; $c_{n-1, -1}$ and $d_{n-1, -1}$ must be interpreted as zero, and the initial values are  $c_{0, 0} = 1$ and $d_{0, 0} = 0$. The system is again obtained by distinguishing the two instances of the last step leading to level $k$.
We set 
\begin{equation*}
\varphi_k=\varphi_k(z)=\sum_{n\ge0}c_{n,k}z^n\quad\text{and}\quad \psi_k=\psi_k(z)=\sum_{n\ge0}d_{n,k}z^n.
\end{equation*}
The translation of the recurrences to generating functions is
\begin{align*}
	\varphi_{ k} & = z\varphi_{ k-1} + z\psi_{ k}+[k=0],\\
	\psi_{ k} & = z\psi_{ k-1} + z\varphi_{ k+1}.
\end{align*}

The matrix of interest is the transposed version of before:
\begin{equation*}
	\mathcal{M}=\begin{pmatrix}
		1&-z^2&&&\\
		-2z&1&-z^2\\
		z^2&	-2z&1&-z^2\\
		&z^2&	-2z&1&-z^2\\
		&&z^2&	-2z&1&-z^2\\
	\end{pmatrix}
\end{equation*}
The equations of interest are
\begin{equation*}
\mathcal{M}\begin{pmatrix}
	\varphi_0\\
	\varphi_1\\
	\varphi_2\\
	\varphi_3\\
\end{pmatrix}=\begin{pmatrix}
	1\\
	-z\\
	0\\
	0\\
\end{pmatrix}\qquad\text{and}\qquad
	\mathcal{M}\begin{pmatrix}
		\psi_0\\
		\psi_1\\
		\psi_2\\
		\psi_3\\
	\end{pmatrix}=\varphi_0\begin{pmatrix}
		1\\
		0\\
		0\\
		0\\
	\end{pmatrix}
\end{equation*}
It is easier to work with the second version, and using that $\varphi_0=\frac{z^2}{1-t}=f_0$, which follows from a combinatorial argument about reading the same combinatorial objects once from left to right resp.\ from right to left. If one wants to avoid such an argument, one can work with the first system, which is similar but a bit messier. We work with $(1,0,0,\dots)^T$ and multiply the solution later by the factor $\frac{z^2}{1-t}$.

We want to apply Cramer's rule to solve the system. Each component will then be a quotient of two determinants. The determinant in the denominator is always the same and the same as before.
Now let $D_{n,i}=D_{n,i}(z)$ be the determinant with $n$ rows, and index $i$ is replaced by a single 1 in the first row and otherwise 0 in row 1 and column $i$. 
These numbers satisfy a recursion, which can be obtained by hand, but it is easier to use a computer (gfun and Maple):
\begin{equation*}
	D_{n,i}=z^{i-1}\tau_{i}D_{n-i}-z^{i+2}\tau_{i-1}D_{n-i-1},
\end{equation*}
where
\begin{equation*}
	\sum_{i\ge0}\tau_{i}w^i=\frac{w}{1-2w+w^2-z^3w^3}.
\end{equation*}
We have this explicit form
\begin{equation*}
\tau _{i}=\sum_{0\le l\le (i-1)/3}\binom{i-l}{2l+1}z^{3l}
\end{equation*}
which is beautiful, but will not be used.

For the application of Cramer's rule we must consider
\begin{equation*}
	\frac{D_{n,i}}{D_n}=z^{i-1}\tau_{i}\frac{D_{n-i}}{D_n}-z^{i+2}\tau_{i-1}\frac{D_{n-i-1}}{D_n}
\end{equation*}
and its limit for $n\to\infty$, which follows from
\begin{equation*}
	\lim_{n\to\infty}\frac{D_{n-i}}{D_n}=(t-1)^{-2i}.
\end{equation*}
Consequently we get
\begin{equation*}
	\lim_{n\to\infty}\frac{D_{n,i}}{D_n}=	z^{i-1}\tau_{i}(t-1)^{-2i}-z^{i+2}\tau_{i-1}(t-1)^{-2i-2}
	=\frac{t^i}{z^{2i+1}}(\tau_{i}-t\tau_{i-1}).
\end{equation*}
But this expression is simpler than $\tau_i$ alone:
\begin{align*}
	\sum_{i\ge0}(\tau_{i}-t\tau_{i-1})w^i=\frac{w(1-tw)}{1-2w+w^2-z^3w^3}=
	\frac{1}{1-2w+w^2+tw-2tw^2+t^2w^2}.
\end{align*}
Since the denominator  is only quadratic (in $w$), there is a Binet formula:
\begin{align*}
\tau_{i}-t\tau_{i-1}=\frac{\mu_3^i-\mu_2^i}{W}
\end{align*}
with
\begin{equation*}
\mu_2=\frac{-t+2-W}{2},\quad\text{and}\quad\mu_3=\frac{-t+2+W}{2},\quad\text{and}\quad W=\sqrt{4t-3t^2}=\mu_3-\mu_2.
\end{equation*}

The final answer is then (notice the shift, since $i$ corresponds to $\psi_{ i-1}$)
\begin{equation*}
\psi_i=\frac{z^2}{1-t}\frac{t^{i+1}}{z^{2i+3}}\frac{\mu_3^{i+1}-\mu_2^{i+1}}{W}=
\frac{t^{i+1}}{z^{2i+1}(1-t)}\frac{\mu_3^{i+1}-\mu_2^{i+1}}{\mu_3-\mu_2}.
\end{equation*}
From
\begin{equation*}
\varphi_{ k}=\frac1z\psi_{ k-1}  - \psi_{ k-2} ,
\end{equation*} 
the values $\varphi_{ 1},\varphi_{ 2},\dots$ can be computed; $\varphi_{ 0}$ is already known, and
consistent with $\varphi_{ 0}=1+z\psi_0=\frac1{1-t}$.

\section{Computing coefficients}

For completeness, we discuss how to compute the coefficients of the generating functions that we obtained. This is more or less a repetition of what we did already in our earlier paper \cite{alone}.

The first ingredient is the inversion of $x = t(1-t)^2$, i. e., $t$ expressed as a series in $x$:
\begin{align*}
[x^n] t^k & = \frac{k}{n}[w^{n-k}]\frac{1}{(1-w)^{2n}} = \frac{k}{n}\binom{3n-k-1}{n-k} \qquad \Rightarrow \qquad t^k = \sum_{n\ge k}\binom{3n-k-1}{n-k}\frac{k}{n}x^n.
\end{align*}
This can be done using the Lagrange inversion formula, or, equivalently by contour integration:
\begin{align*}
[x^n] t^k & = \frac1{2\pi i}\oint\frac{dx}{x^{n+1}}t^k
 = \frac1{2\pi i}\oint\frac{(1-t)(1-3t)dt}{t^{n+1}(1-t)^{2n+2}}t^k\\
&=[t^{n-k}]\frac{1-3t}{(1-t)^{2n+1}}
=\binom{3n-k}{n-k}-3\binom{3n-k-1}{n-k-1}=\frac{k}{n}\binom{3n-k-1}{n-k}.
\end{align*}
The contours are in all instances small circles around the origin.

And now
\begin{align*}
	a_{n,k}=[z^n]f_k(z)&=[z^n]\frac{t^k}{z^k(1-t)}=[z^{n+k}]\frac{t^k}{1-t}=  [x^{(n+k)/3}]\frac{t^k}{1-t}\\
	&=\frac1{2\pi i}\oint\frac{dx}{x^{(n+k)/3+1}}\frac{t^k}{(1-t)}
	=\frac1{2\pi i}\oint\frac{(1-t)(1-3t)dt}{t^{(n+k)/3+1}(1-t)^{2(n+k)/3+2}}\frac{t^k}{(1-t)}\\
	&=[t^{(n-2k)/3}]\frac{1-3t}{(1-t)^{2(n+k)/3+2}}\\
	&=\binom{n+1}{(n-2k)/3-1}-3\binom{n}{(n-2k)/3-2}.
\end{align*}
This works   if $n\equiv 2k \bmod 3$ and $n\ge k$; other values are zero. Note that the forms given in \cite{alone} are different but equivalent.
Similarly 
\begin{align*}
b_{n,k}=[z^n]g_k(z)&=[z^n]\frac{t^{k+1}}{z^{k+2}}=[x^{(n+k+2)/3}]t^{k+1} \\
&=\frac1{2\pi i}\oint\frac{dx}{x^{(n+k+2)/3+1}}t^{k+1}
=\frac1{2\pi i}\oint\frac{(1-t)(1-3t)dt}{t^{(n+k+2)/3+1}(1-t)^{2(n+k+2)/3+2}}t^{k+1}\\
&=[t^{(n-2k-1)/3}]\frac{1-3t}{(1-t)^{2(n+k+2)/3+1}}\\
&=\binom{n+1}{(n-2k-1)/3}-3\binom{n}{(n-2k-1)/3-1}.
\end{align*}
This works   if $n\equiv 2k+1 \bmod 3$ and $n\ge 2k+1$; other values are zero.

Now we show how to extract
\begin{equation*}
d_{n,k}=[z^n]\psi_k(z)=[z^n]\frac{t^{k+1}}{z^{2k+1}(1-t)}\frac{\mu_3^{k+1}-\mu_2^{k+1}}{\mu_3-\mu_2}.
\end{equation*}

The identity \cite[eq. (22)]{Gould} will be useful for the  calculation; it is classical and goes by the name of Girard-Waring
\begin{align*}
	\frac{\mu_3^{k+1}-\mu_2^{k+1}}{\mu_3-\mu_2}
	& = \sum_{i = 0}^{\lfloor k/2 \rfloor}(-1)^i\binom{k-i}{i}(\mu_2+\mu_3)^{k-2i}(\mu_2\mu_3)^i\\*
	& = \sum_{i = 0}^{\lfloor k/2 \rfloor}(-1)^{i+k}\binom{k-i}{i}(t-2)^{k-2i}(t-1)^{2i}.
\end{align*} 
Therefore
\begin{align*}
d_{n,k}=&[z^{n+2k+1}]	 t^{k+1}\sum_{i = 0}^{\lfloor k/2 \rfloor}(-1)^{i+k+1}\binom{k-i}{i}(t-2)^{k-2i}(t-1)^{2i-1}\\
=&[z^{n+2k+1}]	 t^{k+1}\sum_{i = 0}^{\lfloor k/2 \rfloor}(-1)^{i+k+1}\binom{k-i}{i}\sum_{j=0}^{k-2i}\binom{k-2i}{j}(-1)^{k-2i-j}(t-1)^{2i+j-1}\\
=&[z^{n+2k+1}]	 t^{k+1}\sum_{i = 0}^{\lfloor k/2 \rfloor}\sum_{j=0}^{k-2i}(-1)^{i+j+1} \binom{k-i}{i}\binom{k-2i}{j}(t-1)^{2i+j-1}. 
\end{align*} 
To keep this short, we only compute
\begin{align*}
[z^{n+2k+1}]&	 t^{k+1}(t-1)^{M}
=\frac1{2\pi i}\oint\frac{dx}{x^{(n+2k+1)/3+1}}t^{k+1}(t-1)^{M}\\
&=\frac1{2\pi i}\oint\frac{(1-t)(1-3t)dt}{t^{(n+2k+1)/3+1}(1-t)^{2(n+2k+1)/3+2}}t^{k+1}(t-1)^{M}\\
&= (-1)^M[t^{(n-k+1)/3-1}]\frac{1-3t}{(1-t)^{2(n+2k+1)/3-M+1}}\\
&=(-1)^M\bigg[\binom{n+k-M}{(n-k+1)/3-1}-3\binom{n-k-M-1}{(n-k+1)/3-2}\biggr].
\end{align*}
This works  for $n\equiv k-1\bmod 3$. The explicit formul\ae\ for $c_{n,k}$ and
$d_{n,k}$ are given in \cite{alone}.

\begin{acka}
The insightful comments and suggestions of two reviewers are gratefully acknowledged.
	\end{acka}

\bibliographystyle{plain}

\end{document}